\title{A Proof of Green's Conjecture Regarding the Removal Properties of Sets of Linear Equations}
\author{Asaf Shapira
\thanks{Microsoft Research. Email: asafico@tau.ac.il.} }

\date{}

\documentclass [letterpaper,11pt]{article}
\usepackage{amsfonts}
\newtheorem{theo}{Theorem}%[section]

\newtheorem{conj}{Conjecture}

\newtheorem{lemma}{Lemma}[section]
\newtheorem{definition}[lemma]{Definition}

\newtheorem{claim}[lemma]{Claim}

\newcommand{\qed}{\hspace*{\fill} \rule{7pt}{7pt}}

\newcommand{\ignore}[1]{}

\begin{document}
\maketitle

\begin{abstract}

A system of $\ell$ linear equations in $p$ unknowns $Mx=b$ is said
to have the {\em removal property} if every set $S \subseteq
\{1,\ldots,n\}$ which contains $o(n^{p-\ell})$ solutions of $Mx=b$
can be turned into a set $S'$ containing no solution of $Mx=b$, by
the removal of $o(n)$ elements. Green [GAFA 2005] proved that a
single homogenous linear equation always has the removal property,
and conjectured that every set of homogenous linear equations has
the removal property. We confirm Green's conjecture by showing that
every set of linear equations (even non-homogenous) has the removal
property.

\end{abstract}

\section{Introduction}\label{intro}

The (triangle) removal lemma of Ruzsa and Szemer\'edi \cite{RuS},
which is by now a cornerstone result in combinatorics, states that a
graph on $n$ vertices that contains only $o(n^3)$ triangles can be
made triangle free by the removal of only $o(n^2)$ edges. Or in
other words, if a graph has asymptomatically few triangles then it
is asymptotically close to being triangle free. While the lemma was
proved in \cite{RuS} for triangles, an analogous result for any
fixed graph can be obtained using the same proof idea. Actually, the
main tool for obtaining the removal lemma is Szemer\'edi's
regularity lemma for graphs \cite{Sz}, another landmark result in
combinatorics. The removal lemma has many applications in different
areas like extremal graph theory, additive number theory and
theoretical computer science. Perhaps its most well known
application appears already in \cite{RuS} where it is shown that an
ingenious application of it gives a very short and elegant proof of
Roth's Theorem, which states that every $S \subseteq
[n]=\{1,\ldots,n\}$ of positive density contains a 3-term arithmetic
progression.

Recall that an $r$-uniform hypergraph $H=(V,E)$ has a set of
vertices $V$ and a set of edges $E$, where each edge $e \in E$
contains $r$ distinct vertices from $V$. So a graph is a 2-uniform
hypergraph. Szemeredi's famous theorem \cite{Sztheo} extends Roth's
theorem by showing that every $S \subseteq [n]$ of positive density
actually contains arbitrarily long arithmetic progressions (when $n$
is large enough). Motivated by the fact the a removal lemma for
graphs can be used to prove Roth's theorem, Frankl and R\"odl
\cite{FR} showed that a removal lemma for $r$-uniform hypergraphs
could be used to prove Szemeredi's theorem on $(r+1)$-term
arithmetic progressions. They further developed a regularity lemma,
as well as a corresponding removal lemma, for 3-uniform hypergraphs
thus obtaining a new proof of Szemeredi's theorem for 4-term
arithmetic progressions. In recent years there have been many
exciting results in this area, in particular the results of Gowers
\cite{Gowers} and of Nagle, R\"odl Schacht and Skokan \cite{NRS,RS},
who independently obtained regularity lemmas and removal lemmas for
$r$-uniform hypergraph, thus providing alternative combinatorial
proofs of Szemeredi's Theorem \cite{Sztheo} and some of it
generalizations, notably those of Furstenberg and Katznelson
\cite{FK}. Tao \cite{Tao} later obtained another proof of the
hypergraph removal lemma and of its many corollaries mentioned
above. For more details see \cite{Gowers2,KNRSS}.

In this paper we will use the above mentioned hypergraph removal
lemma in order to resolve a conjecture of Green \cite{Green}
regarding the removal properties of sets of linear equations. Let
$Mx=b$ be a set of linear equations, and let us say that a set of
integers $S$ is $(M,b)$-free if it contains no solution to $Mx=b$,
that is, if there is no vector $x$, whose entries all belong to $S$,
which satisfies $Mx=b$. Just like the removal lemma for graphs
states that a graph that has few copies of $H$ should be close to
being $H$-free, a removal lemma for sets of linear equations $Mx=b$
should say that a subset of the integers $[n]$ that contains few
solutions to $Mx=b$, should be close to being $(M,b)$-free. Let us
start be defining this notion precisely.

\begin{definition}[Removal Property]\label{removalproperty} Let $M$
be an $\ell \times p$ matrix of integers and let $b \in
\mathbb{N}^{\ell}$. The set of linear equations $Mx=b$ has the {\em
removal property} if for every $\delta >0$ there is an
$\epsilon=\epsilon(\delta,M,b)>0$ with the following property: if $S
\subseteq [n]$ is such that there are at most $\epsilon n^{p-\ell}$
vectors $x \in S^{p}$ satisfying $Mx=b$, then one can remove from
$S$ at most $\delta n$ elements to obtain an $(M,b)$-free set.
\end{definition}

We note that in the above definition, as well as throughput the
paper, we assume that the $\ell \times p$ matrix $M$ of
a set of linear equations has rank $\ell$.

Green \cite{Green} has initiated the study of the removal properties
of sets of linear equations. His main result was the following:

\begin{theo}[Green \cite{Green}]\label{theogreen} Any single homogenous
linear equation has the removal property.
\end{theo}

The main result of Green actually holds over any abelian group. To
prove this result, Green developed a regularity lemma for abelian
groups, which is somewhat analogous to Szemer\'edi's regularity
lemma for graphs \cite{Sz}. Although the application of the group
regularity lemma for proving Theorem \ref{theogreen} was similar to
the derivation of the graph removal lemma from the graph regularity
lemma, the proof of the group regularity lemma was far from trivial.
One of the main conjectures raised in \cite{Green} is that a natural
generalization of Theorem \ref{theogreen} should also hold
(Conjecture 9.4 in \cite{Green}).

\begin{conj}[Green \cite{Green}]\label{Greensconj} Any system of
homogenous linear equations $Mx=0$ has the removal
property.
\end{conj}

We note that besides being a natural generalization of Theorem
\ref{theogreen}, Conjecture \ref{Greensconj} was also raised in
\cite{Green} with relation to a conjecture of Bergelson, Host, Kra
and Ruzsa \cite{BHK} regarding the number of $k$-term arithmetic
progressions with a common difference in subsets of $[n]$. See
Section \ref{concluding} for more details.

Very recently, Kr\'al', Serra and Vena \cite{Kral} gave a
surprisingly simple proof of Theorem \ref{theogreen}, which
completely avoided the use of Green's regularity lemma for groups.
In fact, their proof is an elegant and simple application of the
graph removal lemma mentioned earlier and it actually extends
Theorem \ref{theogreen} to any single non-homogenous linear equation
over non-abelian groups. Kr\'al', Serra and Vena \cite{Kral} also
show that Conjecture \ref{Greensconj} holds when $M$ is a 0/1
matrix, which satisfies certain conditions. But these conditions are
not satisfied even by all 0/1 matrices. In another recent result, 
which was obtained independently of ours, Candela \cite{C} showed
that Conjecture \ref{Greensconj} holds for every pair of homogenous 
linear equations, as well as for every system of homogenous equations in 
which every $\ell$ columns of $M$ are linearly independent. 
See more details in Subsection \ref{subsecover}.

In this paper we confirm Green's for every homogenous set of linear
equations. In fact, we prove the following more general result.

\begin{theo}[Main Result]\label{maintheo} Any set of linear equations
(even non homogenous) $Mx=b$ has the removal property.
\end{theo}

The rest of the paper if organized as follows. In the next section
we give an overview of the proof of Theorem \ref{maintheo}. As we
show in that section, Theorem \ref{maintheo} also holds over any
finite field, that is when $S \subseteq \mathbb{F}_n$, where
$\mathbb{F}_n$ is the field of size $n$. In fact it is easy to
modify the proof so that it works over any field, but we will not do
so here. The proof of Theorem \ref{maintheo} has two main steps: the first one,
described in Lemma \ref{removal4linear}, applies the main idea from
\cite{Kral} in order to show that if a set of linear equations can
be ``represented'' by a hypergraph then Theorem \ref{maintheo} would
follow from the hypergraph removal lemma. So the second, and most
challenging step of the proof, is showing that every set of linear
equations can be represented as a hypergraph. The proof of this
result, stated in Lemma \ref{mainstep}, appears in Section
\ref{secmain}. In Section \ref{concluding} we give some concluding
remarks and discuss some open problems.

\section{Proof Overview}\label{SecEquiv}

It will be more convenient to deduce Theorem \ref{maintheo} from an
analogous result over the finite field $\mathbb{F}_n$ of size $n$
(for $n$ a prime power). In fact, somewhat surprisingly, we will
actually need to prove a stronger claim than the one asserted in
Theorem \ref{maintheo}. This more general variant, stated in Theorem
\ref{removaloverff}, allows each of the variables $x_i$ to have its
own subset $S_i \subseteq [n]$. We note that a proof of this variant
of Theorem \ref{maintheo} for the case of a single equation was
already proved in \cite{Green} and \cite{Kral}, but in those papers
it was not necessary to go through this more general result. As we
will explain later (see Claims \ref{notwo} and \ref{hsimple}), the
fact that we are considering a more general problem will allow us to
overcome some degeneracies in the system of equations by allowing us
to remove certain equations. This manipulation can be performed when
one considers the generalized removal property (defined below) but
there is no natural way of performing these manipulations when
considering the standard removal property. Therefore, proving this
extended result is essential for our proof strategy.

In what follows and throughout the paper, whenever $x$ is a vector,
$x_i$ will denote its $i^{th}$ entry. Similarly, if $x_1,\ldots,x_p$
are elements in a field, then $x$ will be the vector whose entries
are $x_1,\ldots,x_p$. We say that a collection of $p$ subsets
$S_1,\ldots,S_p \subseteq \mathbb{F}_n$ is $(M,b)$-free if there are
no $x_1 \in S_1,\ldots, x_p \in S_p$ which satisfy $Mx=b$.

\begin{definition}[Generalized Removal Property over Finite
Fields]\label{removalpropertyFF} Let $\mathbb{F}_n$ be the field of
size $n$, let $M$ be an $\ell \times p$ matrix over $\mathbb{F}_n$
and let $b \in \mathbb{F}^{\ell}_n$. The system $Mx=b$ is said to
have the {\em generalized removal property} if for every $\delta >0$
there is an $\epsilon=\epsilon(\delta,p)>0$ such that if
$S_1,\ldots,S_p \subseteq \mathbb{F}_n$ contain less than $\epsilon
n^{p-\ell}$ solutions to $Mx=b$ with each $x_i \in S_i$, then one
can remove from each $S_i$ at most $\delta n$ elements to obtain
sets $S'_1,\ldots,S'_p$ which are $(M,b)$-free.
\end{definition}

By taking all sets $S_i$ to be the same set $S$ we, of course, get
the standard notion of the removal property from Definition
\ref{removalproperty} so we may indeed work with this generalized
definition. We will deduce Theorem \ref{maintheo} from the following
theorem.

\begin{theo}\label{removaloverff}
Every set of linear equations
$Mx=b$ over a finite field has the generalized removal property.
\end{theo}

In this paper we apply the hypergraph removal lemma in order to
resolve Green's conjecture. In fact, for the proof of Theorem
\ref{removaloverff} we will need a variant of the hypergraph removal
lemma which works for colored hypergraphs. But let us first recall
some basic definitions. An $r$-uniform hypergraph is simple if it
has no parallel edges, that is, if different edges contain different
subsets of vertices of size $r$. We say that a set of vertices $U$
in a $r$-uniform hypergraph $H=(V_H,E_H)$ span a copy of an
$r$-uniform hypergraph $K=(V_K,E_K)$ if there is an injective
mapping $\phi$ from $V_K$ to $U$ such that if $v_1,\ldots,v_r$ form
an edge in $K$ then $\phi(v_1),\ldots,\phi(v_r)$ form an edge in $U
\subseteq V_H$. We say that a hypergraph is $c$-colored if its edges
are colored by $\{1,\ldots,c\}$. If $K$ and $H$ are $c$-colored,
then $U$ is said to span a colored copy of $K$ if the above mapping
$\phi$ sends edges of $K$ of color $i$ to edges of $H$ (in $U$) of
the same color $i$. We stress that the coloring of the edges does
not have to satisfy any constraints that are usually associated with
edge colorings. Finally, the number of colored copies of $K$ in $H$
is the number of subsets $U \subseteq V_H$ of size $|V_K|$ which
span a colored copy of $K$.

The following variant of the hypergraph removal lemma is a special
case of Theorem 1.2 in \cite{AT}.\footnote{As noted to us by Terry
Tao, this variant of the hypergraph removal lemma can probably be
extracted from the previous proofs of the hypergraph removal lemma
\cite{Gowers,NRS,RS,Tao}, just like the colored removal lemma for
graphs can be extracted from the proof of the graph removal lemma,
see \cite{KS}.}

\begin{theo}[Austin and Tao \cite{AT}]\label{Cremoval} Let $K$ be a fixed
$r$-uniform $c$-colored hypergraph on $k$ vertices. For every
$\delta >0$ there is an $\epsilon=\epsilon(\delta,k)>0$ such that if
$H$ is an $r$-uniform $c$-colored simple hypergraph with less than
$\epsilon n^k$ colored copies of $K$, then one can remove from $H$
at most $\delta n^r$ edges and obtain a hypergraph that contains no
colored copy of $K$.
\end{theo}

In order to use Theorem \ref{Cremoval} for the proof of Theorem
\ref{removaloverff}, we will need to represent the solutions of
$Mx=b$ as colored copies of a certain ``small'' hypergraph $K$ in a
certain ``large'' hypergraph $H$. The following notion of hypergraph
representability specifies the requirements from such a
representation that suffice for allowing us to deduce Theorem
\ref{removaloverff} from Theorem \ref{Cremoval}.

\begin{definition}[Hypergraph Representation]\label{DefHyperRep}
Let $\mathbb{F}_n$ be the field of size $n$, let $M$ be an $\ell
\times p$ matrix over $\mathbb{F}_n$. The system of linear equations
$Mx=b$ is said to be {\em hypergraph representable} if there is an
integer $r=r(M,b) \leq p^2$ and an $r$-uniform $p$-colored
hypergraph $K$ with $k=r-1+p-\ell$ vertices and $p$ edges, such that
for any $S_1,\ldots,S_p \subseteq [n]$ there is an $r$-uniform
hypergraph $H$ on $k n$ vertices which satisfies the following:
\begin{enumerate}
\item $H$ is simple and each edge with color $i$ is labeled by one of the elements of
$S_i$.

\item If $x_1 \in S_1,\ldots,x_p \in S_p$ satisfy $Mx=b$
then $H$ contains $n^{r-1}$ colored copies of $K$, such that their
edge with color $i$ has label $x_i$. These colored copies of $K$
should also be edge disjoint.

\item If $S_1,\ldots,S_p$ contain $T$ solutions to $Mx=b$ with
$x_i \in S_i$ then $H$ contains $Tn^{r-1}$ colored copies of $K$.

\end{enumerate}
\end{definition}

The following lemma shows that a hypergraph representation can allow
us to prove Theorem \ref{removaloverff} using the hypergraph removal
lemma.

\begin{lemma}\label{removal4linear} If $Mx=b$ has a hypergraph
representation then it has the generalized removal property.
\end{lemma}

\paragraph{Proof:} Suppose $Mx=b$ is a system of $\ell$ linear equations in
$p$ unknowns. Let $S_1,\ldots,S_p$ be $p$ subsets of $\mathbb{F}_n$
and let $H$ be the hypergraph guaranteed by Definition
\ref{DefHyperRep}. We claim that we can take $\epsilon(\delta,p)$ in
Theorem \ref{removalpropertyFF} to be the value
$\epsilon=\epsilon(\delta/pk^{r},k)$ from Lemma \ref{Cremoval}. Note
that $r,k \leq 2p^2$ so this still implies that $\epsilon$ is only a
function of $\delta$ and $p$. Indeed, if $S_1,\ldots,S_p$ contain
only $\epsilon n^{p-\ell}$ solutions to $Mx=b$ then by item 3 of
Definition \ref{DefHyperRep} we get that $H$ contains at most
$\epsilon n^{p-\ell} \cdot n^{r-1}=\epsilon n^k$ colored copies of
$K$. As $H$ is simple, we can apply the removal lemma for colored
hypergraphs (Lemma \ref{Cremoval}) to conclude that one can remove a
set $E$ of at most $\frac{\delta}{pk^r} (kn)^r=\frac{\delta}{p} n^r$
edges from $H$ and thus destroy all the colored copies of $K$ in $H$
(recall that $H$ has $kn$ vertices).

To show that we can turn $S_1,\ldots,S_p$ into a collection of
$(M,b)$-free sets by removing only $\delta n$ elements from each
$S_i$, let us remove an element $s$ from $S_i$ if $E$ contains at
least $n^{r-1}/p$ edges that are colored with $i$ and labeled with
$s$. As each edge has one label (because $H$ has no parallel edges),
and $|E| \leq \frac{\delta}{p} n^r$ this means that we remove only
$\delta n$ elements from each $S_i$. To see that we thus turn
$S_1,\ldots,S_p$ into $(M,b)$-free sets, suppose that the new sets
$S'_1,\ldots,S'_p$ still contain a solution $s_1 \in S_1,\ldots,s_p
\in S_p$ to $Mx=b$. By item 2 of Definition \ref{DefHyperRep}, this
solution defines $n^{r-1}$ edge disjoint colored copies of $K$ in
$H$, with the property that in every colored copy, the edge with
color $i$ is labeled with the same element $s_i \in S_i$. As $E$
must contain at least one edge from each of these colored copies (as
it should destroy all such copies), there must be some $1 \leq i
\leq p$ for which $E$ contains at least $n^{r-1}/p$ edges that are
colored $i$ and labeled with $s_i$. But this contradicts the fact
that $s_i$ should have been removed from $S_i$. $\qed$

\bigskip

We note that the above lemma generalizes a similar lemma for the
case of representing a single equation using a graph, which was
implicit in \cite{Kral}. In fact, as we have mentioned earlier,
\cite{Kral} also show that a set of homogenous linear equations
$Mx=0$, with $M$ being a 0/1 matrix, that satisfies certain
conditions also has the removal lemma. One of these conditions
essentially says that the system of equations is {\em graph}
representable. However, there are even some 0/1 matrices for which
$Mx=0$ is not graph representable (in the sense of \cite{Kral}).
Lemma \ref{mainstep} below shows that any set of linear equations
has a {\em hypergraph} representation. This lemma is proved in the
next section and it is the most challenging part of this paper.

\begin{lemma}\label{mainstep} Every set of linear equations $Mx=b$
over a finite field is hypergraph representable.
\end{lemma}

From the above two lemmas we get the following.

\paragraph{Proof of Theorem \ref{maintheo}:} Immediate from
Theorem \ref{removaloverff} and Lemma \ref{removal4linear}.

\bigskip

As we have mentioned before, Theorem \ref{mainstep} is now an easy
application of Theorem \ref{removaloverff}.

\paragraph{Proof of Theorem \ref{maintheo}:} Given a set of
linear equations $Mx=b$ in $p$ unknowns, let $c$ be the maximum
absolute value of the entries of $M$ and $b$. Given an integer $n$
let $q=q(n)$ be the smallest prime larger than $cp^2n$. It is well
known that $q \leq 2cp^2n$ (in fact, much better bounds are known).
It is clear that for a vector $x \in [n]^p$ we have $Mx=b$ over
$\mathbb{R}$ if and only if $Mx=b$ over $\mathbb{F}_q$. So if $Mx=b$
has $o(n^{p-\ell})$ solutions with $x_i \in S_i$ over $\mathbb{R}$,
it also has $o(q^{p-\ell})$ solutions with $x_i \in S_i \subseteq
\mathbb{F}_q$ over $\mathbb{F}_q$. By Theorem \ref{removaloverff} we
can remove $o(q)$ elements from each $S_i$ and obtain sets $S'_i$
that are $(M,b)$-free. But as $q=O(n)$ we infer that the removal of
the same $o(q)=o(n)$ elements also guarantees that the sets are
$(M,b)$-free over $\mathbb{R}$. $\qed$

\subsection{Overview of the Proof of Lemma \ref{mainstep}}\label{subsecover}

Let us start by noting that Lemma \ref{mainstep} for the case of a
single equation was (implicity) proven in \cite{Kral}, where they
show that one can take $r=2$, in other words, they represent a
single equation as a graph $K$, in a graph $H$. Actually, the graph
$K$ in the proof of \cite{Kral} is a cycle of length $p$. The proof
in \cite{Kral} is very short and elegant, and we recommend reading
it to better understand the intuition behind our proof (although
this paper is, of course, self contained). Another related result is
the proof of Szemer\'edi's theorem \cite{Sztheo} using the
hypergraph removal lemma \cite{FR}, which can be interpreted as
(essentially) showing that the set of $p-2$ linear equations which
define a $p$-term arithmetic progression\footnote{These linear
equations are $x_1+x_3=2x_2$,
$x_2+x_4=2x_3,\ldots,x_{p-2}+x_p=2x_{p-1}$.} are hypergraph
representable with $K$ being the complete $(p-1)$-uniform hypergraph
of size $p$. ``Interpolating'' these two special cases of Lemma \ref{mainstep}
suggests that a hypergraph representation of a set of $\ell$ linear
equations in $p$ unknowns should involve an $(\ell+1)$-uniform
hypergraph $K$ of size $p$. And indeed, we initially found a
(relatively) simple way to achieve this for $p-2$ equations in $p$
unknowns, thus extending the representability of the arithmetic
progression set of linear equations.

However, somewhat surprisingly, when $1 < \ell < p-2$ the situation
becomes much more complicated and we did not manage to find a simple
representation along the lines of the above two cases. The problem with trying
to extend the previous approaches to larger sets of equations is that
obtaining all the requirements of Definition \ref{DefHyperRep} turns
out to be very complicated when $M$ has a set of $\ell$ columns that 
are not linearly independent. Let us mention again that Candela \cite{C}
has recently considered linear equations $Mx=0$ in which every $\ell$
columns are linearly independent, and showed that Conjecture \ref{Greensconj} 
holds in these cases.

The way we overcome the above complications is by using a representation involving
hypergraphs of a much larger degree of uniformity (that is, larger
edges), which is roughly the number of non-zero entries of $M$ after
we perform certain manipulations on it. We note that specializing
our proof to either the case $\ell=1$ or to the case $\ell=p-2$ does
{\em not} give proofs that are identical to the ones (implicit) in
\cite{FR} or \cite{Kral}. For example, our proof for the case of a
single equation in $p$ unknowns uses a $(p-1)$-uniform hypergraph,
rather than a graph as in \cite{Kral}.

So let us give a brief overview of the proof. We need to find a
small hypergraph $K$ with $p$ edges, whose copies, within another
hypergraph $H$, will represent the solutions to $Mx=b$. Each edge of
$H$, and therefore also $K$, will have a color $1 \leq i \leq p$ and
a label $s \in S_i$. The system $Mx=b$ has $p$ unknowns and $K$ has
$p$ edges and it may certainly be the case that all the entries of
$M$ are non-zero. It is apparent that using all the edges of $K$ to
``deduce'' a linear equation of $Mx=b$ is not a good idea because in
that way we will only be able to extract one equation from a copy of
$K$ and we need to extract $\ell$ such equations. Therefore, we will
first ``diagonalize'' an $\ell \times \ell$ sub-matrix of $M$ to get
an equivalent set of equations (which we still denote by $Mx=b$)
which has the property that $p-\ell$ of its unknowns
$x_1,\ldots,x_{p-\ell}$ (can) appear in all equations and the rest
of the $\ell$ unknowns $x_{p-\ell+1},\ldots,x_p$ each appear in
precisely one equation. This suggests the idea of extracting
equation $i$ from (some of) the edges corresponding to
$x_1,\ldots,x_{p-\ell}$ and one of the edges corresponding to
$x_{p-\ell+1},\ldots,x_p$. The hypergraph $K$ first contains
$p-\ell$ edges that do not depend on the structure of $M$. The other
$\ell$ edges do depend on the structure of $M$ and use the previous
$p-\ell$ edges in order to ``construct'' the equations of $Mx=b$.
The way to think about this is that for any copy of $K$ in $H$ the
first $p-\ell$ edges will have a special vertex that will hold a
value from $S_i$ (this will be the vertex in one of the sets
$U_1,\ldots,U_{p-\ell}$ defined in Section \ref{secmain}). The other
$\ell$ edges will include some of these special vertices, depending
on the equation we are trying to build. The way we will deduce an
equation from a copy of $K$ in $H$ is that we will argue that the
fact that two edges have a common vertex means that a certain
equation holds. See Claim \ref{edgeequation}.

But there is another complication here because the linear equation
we obtain in the above process will contain many other variables not
from the sets $S_i$, which will need to vanish from such an
equation, in order to allow us to extract the linear equations we
are really interested in. The reason for these ``extra'' variables
is that $H$ needs to contain $n^{r-1}$ edge disjoint copies of $K$
for every solution of $Mx=b$. Hence, an edge of $H$ will actually be
parameterized by several other elements from $\mathbb{F}_n$ (these
are the elements $x_1,\ldots,x_{r-1}$ that are used after Claim
\ref{Bi}). So we will need to make sure that these extra variables
vanish in the linear equation which we extract from a copy of $K$.
To make sure this happens we will need to carefully choose the
vertices of each edge within $H$.

A final complication arises from the fact that while we need $H$ to
contain relatively few copies of $H$, we also need it to contain
many copies edge disjoint copies of $H$ for every solution of
$Mx=b$. To this end we will think of each vertex of $H$ as a linear
equation and we will want the linear equations corresponding to the
vertices of an edge to be linearly independent. The reason why it is
hard to prove Lemma \ref{mainstep} using an $(\ell+1)$-uniform
hypergraph (as the results of \cite{Kral} and \cite{FR} may suggest)
is that it seems very hard to obtain all the above requirements
simultaneously. The fact that we are considering hypergraphs with a
larger degree of uniformity will allow us (in some sense) to break
the dependencies between these requirements.

\section{Proof of Lemma \ref{mainstep}}\label{secmain}

Let $M$ be an $\ell \times p$ matrix over $\mathbb{F}_n$ and $b \in
\mathbb{F}^{\ell}_n$. We will first perform a series of operations
on $M$ and $b$ which will help us in proving Lemma \ref{mainstep}.
For convenience, we will continue to refer to the transformed matrix
and vector as $M$ and $b$. Suppose, without loss of generality, that
the last $\ell$ columns of $M$ are linearly independent. We can thus
transform $M$ (and accordingly also $b$) into an equivalent set of
equations in which the last $\ell$ columns form an identity matrix.
For a row $M_i$ of $M$ let $m_i$ be the largest index $1 \leq j \leq
p-\ell$ for which $M_i$ is non-zero. Let $W_i$ denote the set of
indices $1 \leq j \leq m_i-1$ for which $M_{i,j}$ is non-zero.
Therefore, $M_i$ has $|W_i|+2$ non-zero entries. We will need the
following claim, in which we make use of the fact that we are
actually proving that every set of equations has the generalized
removal property and not just the removal property.

\begin{claim}\label{notwo}
Suppose that every set of $\ell-1$ equations in $p-1$ unknowns
over $\mathbb{F}_n$ has the generalized removal property. Suppose
that the matrix $M$ defined above has a row with less than 3 non-zero
entries. Then $Mx=b$ has the generalized removal property as
well.
\end{claim}

\paragraph{Proof:} Suppose that (say) the first row of $M$ has at
most 2 non-zero entries. If this row has two non-zero elements then
we can assume without loss of generality that it is of the form
$x_1=b-a \cdot x_j$ where $p-\ell+1 \leq j \leq p$. But then we can
get an equivalent set of linear equations $M'x=b'$ by removing the
first row from $M$, removing the column in which $x_j$ appears
(because $x_j$ does not appear in other rows), removing the first
entry of $b$ and updating $S_1$ to be $S'_1=S_1 \cap \{b-a \cdot s
~:~ s \in S_j\}$. We thus get an instance $M'x=b'$ with $\ell-1$
equations and $p-1$ unknowns, hence we can use the assumption of the
claim because: (i) The number of solutions of $Mx=b$ with $x_i \in
S_i$ is precisely the number of solutions of $M'x=b'$ with $x_1 \in
S'_1,x_2\in S_2,\ldots,x_{j-1} \in S_{j-1},x_{j+1} \in
S_{j+1},\ldots,x_p\in S_p$ (ii) if we can remove $\delta n$ elements
from each of the sets of the new instance and thus obtain sets with
no solution of $M'x=b'$ then the removal of the same elements from
the original sets $S_i$ would also give sets with no solution of
$Mx=b$.

If the first row of $M$ has just one non-zero entry, then this
equation is of the form $x_j=b$ for some $p-\ell+1 \leq j \leq p$
and $b \in \mathbb{F}_n$. If $b \notin S_j$ then the sets contain no
solution to $Mx=b$ and there is nothing to prove. If $b \in S_j$
then the number of solutions to $Mx=b$ is the number of solutions of
the set of equations $M'x=b'$ where $M'$ is obtained by removing the
row and column to which $x_j$ belongs and by removing the first
entry of $b$. As in the previous case we can now use the assumption
of the claim. $\qed$

\bigskip

Claim \ref{notwo} implies that we can assume without loss of
generality that none of the sets $W_1,\ldots,W_{\ell}$ is empty,
because if one of them is empty then the corresponding row of $M$
contains less than 3 non-zero entries. In that case we can iteratively remove
equations from $M$ until we either: (i) get a set of linear equations in
which none of the rows has less than 3 non-zero entries, in which case
we can use the fact that the result holds for such sets of equations as we
will next show, or (ii) we get a single equation with only 2 unknowns
with a non-zero coefficient\footnote{Note
that this process can result in having unknowns with a zero coefficient in all
the remaining equations.}. It is now
easy to see that such an equation has the removal property. Indeed,
suppose the equation has $p$ unknowns and only $x_1$ and $x_2$ have a
non-zero coefficient. So the equation is $a_1\cdot x_1 + a_2 \cdot x_2 +\sum^p_{i=3}0\cdot x_i=b$.
In this case the number of solutions to the equation
from sets $S_1,\ldots,S_p$ is the number of solutions to the equation
$a_1x_1+a_2x_2=b$ with $x_1 \in S_1, x_2 \in S_2$ multiplied by $\prod^p_{i=3} |S_i|$.
Therefore, if $S_1,\ldots,S_p$ contain $o(n^{p-1})$ solutions, then either (i) one of the sets
$S_3,\ldots,S_p$ is of size $o(n)$, so we can remove all the elements from this set, or (ii) $S_1,S_2$
contain $o(n)$ solutions to $a_1\cdot x_1+a_2\cdot x_2=b$, but in this case, for every solution $(s_1,s_2)$
we can remove $s_1$ from $S_1$. In either case the new sets $S'_1,\ldots,S'_p$ contain
no solution of the equation, as needed.

We now return to the proof of Lemma \ref{mainstep}, with the assumption
that none of the sets $W_i$ is empty.
Let us multiply each of the rows of $M$ by $M^{-1}_{i,m_i}$ so that
for every $1 \leq i \leq \ell$ we have $M_{i,m_i}=1$. For every $1
\leq i \leq \ell$ let $d_i \in \{p-\ell+1,\ldots, p\}$ denote the
index of the unique non-zero entry of $M_i$ within the last $\ell$
columns of $M$. Using the notation which we have introduced thus
far, the system of linear equations $Mx=b$ can be written as the set
of $\ell$ equations $L_1,\ldots,L_{\ell}$, where $L_i$ is the
equation
\begin{equation}\label{rewriteeq}
x_{m_i}+M_{i,d_i}\cdot x_{d_i} +\sum_{j \in W_i}M_{i,j}\cdot
s_{j}=b_i\;.
\end{equation}

Let us set
$$
r=1+\sum_{1 \leq i \leq \ell}|W_i|\;.
$$
Observe that as mentioned in the statement of the lemma, we indeed
have $r \leq p^2$.

We now define an $r$-uniform $p$-colored hypergraph $K$, which will
help us in proving that $Mx=b$ is hypergraph representable as in
Definition \ref{DefHyperRep}. The hypergraph $K$ has $k=r-1+p-\ell$
vertices which we denote by
$v_1,\ldots,v_{r-1},u_1,\ldots,u_{p-\ell}$. As for $K$'s edges, it
first contains $p-\ell$ edges denoted $e_1,\ldots,e_{p-\ell}$, where
$e_i$ contains the vertices $v_1,\ldots,v_{r-1},u_i$. Note that
these edges do not depend on the system $Mx=b$. As we will see
later, these edges will help us to ``build'' the actual
representation of the linear equations of $Mx=b$. So in addition to
the above $p-\ell$ edges, $K$ also contains $\ell$ edges
$f_{p-\ell+1},\ldots,f_{p}$, where edge $f_{d_i}$ will\footnote{Note
that we are using the fact that $d_1,\ldots,d_{\ell}$ are distinct
numbers in $\{p-\ell+1,\ldots,p\}$.} represent (in some sense)
equation $L_i$, defined in (\ref{rewriteeq}). To define these $\ell$
edges it will be convenient to partition the set $[r-1]$ into $\ell$
subsets $I_1,\ldots,I_{\ell}$ such that $I_1$ contains the numbers
$1,\ldots,|W_1|$, and $I_2$ contains the numbers
$|W_1|+1,\ldots,|W_1|+|W_2|$ and so on. With this partition we
define for every $1 \leq i \leq \ell$ edge $f_{d_i}$ to contain the
vertices $\{v_i ~:~ i \in [r-1] \setminus I_i\}$, the vertices
$\{u_j~:~ j \in W_i\}$ as well as vertex $u_{m_i}$. Note that as
$|I_i|=|W_i|$ the hypergraph $K$ is indeed $r$-uniform. As for the
coloring of the edges of $K$, for every $1 \leq i \leq p-\ell$ edge
$e_i$ is colored $i$ and for every $p-\ell+1 \leq d_i \leq p$ edge
$f_{d_i}$ is colored $d_i$.

Before defining the hypergraph $H$ we need to define $p-\ell$
vectors $a^1,\ldots,a^{p-\ell} \in \mathbb{F}^{r-1}_n$ which we will
use when defining $H$. We think of $a^1,\ldots,a^{p-\ell}$ as the
$p-\ell$ rows of a $p-\ell \times r-1$ matrix $A$. Furthermore, for
every $1 \leq i \leq p-\ell$ let $A_i$ be the sub-matrix of $A$
which contains the columns whose indices belong to $I_i$ (which was
defined above). We now take the (square) sub-matrix of $A_i$ which
contains the rows whose indices belong to $W_i$ to be the identity
matrix (over $\mathbb{F}_n$). More precisely, if the elements of
$W_i$ are $j_1 < j_2 < \ldots < j_{|W_i|}$ then $A'_{j_g,g}=1$ for
every $1 \leq g \leq |W_i|$, and 0 otherwise\footnote{Note that the
second index of $A'_{j_g,g}$ refers to the column number within
$A_i$, not $A$.}. For future reference, let's denote by $A'_i$ this
square sub-matrix of $A_i$. We finally set row $m_i$ of $A_i$ to be
the vector whose $g^{th}$ entry is $-M_{i,j_g}$, where as above
$j_g$ is the $g^{th}$ element of $W_i$. If $A_i$ has any other rows
besides the ones defined above, we set them to $0$. As each column
of $A$ belongs to one of the matrices $A_i$ we have thus defined $A$
and therefore also the vectors $a^1,\ldots,a^{p-\ell}$.

Let us make two simple observations regarding the above defined
vectors which we will use later. First, let $1 \leq i \leq \ell$ and
$t \in I_i$ and suppose $t$ is the $g^{th}$ element of $I_i$.
Then\footnote{Note that $t$ is an index of a column of $A$, while
$g$ is an index of a column of $A_i$.}
\begin{equation}\label{nicesum}
\sum_{j \in W_i}a^{j}_t \cdot M_{i,j}=(A_i)_{j_g,g} \cdot M_{i,j_g}=M_{i,j_g}=-(A_i)_{m_i,g}=-a^{m_i}_t\;,
\end{equation}
where the first equality is due to the fact that the only non-zero
entries within column $g$ of $A_i$ and the rows from $W_i$ appears
in row $j_g$. The second equality uses the fact that this entry is
in fact 1. The third equality uses the definition of row $m_i$ of
$A_i$.

The second observation we will need is the following.

\begin{claim}\label{Bi} For $1 \leq i \leq \ell$, let $B_i$ be the
following $r-1 \times r-1$ matrix: for every $j \in [r-1] \setminus
I_i$ we have $(B_i)_{j,j}=1$ and $(B_i)_{j,t}=1$ for $t \neq j$. The other
$|I_i|$ rows of $B_i$ are the $|W_i| ~(=|I_i|)$ vectors $\{a^t : t \in
W_i\}$. Then, for every $1 \leq i \leq \ell$ the matrix $B_i$ is
non-singular.
\end{claim}

\paragraph{Proof:} To show that $B_i$ is non-singular it is clearly
enough to show that its $|I_i| \times |I_i|$ minor $B'_i$, which is
determined by $I_i$, is non-singular. But observe that this fact
follows from the way we have defined the vectors
$a^1,\ldots,a^{p-\ell}$ above because $B'_i$ is just $A'_i$, which
is in fact the identity matrix. $\qed$

\bigskip

We are now ready to define, for every set of subsets $S_1,\ldots,S_p
\subseteq \mathbb{F}_n$, the hypergraph $H$ which will establish
that $Mx=b$ is hypergraph representable. The vertex set of $H$
consists of $k~(=r-1+p-\ell)$ disjoint sets
$V_1,\ldots,V_{r-1},U_1,\ldots,U_{p-\ell}$, where each of these sets
contains $n$ vertices and we think of the elements of each of these
sets as the elements of $\mathbb{F}_n$. As for the edges of $H$, we
first put for $1 \leq i \leq p-\ell$ and every choice of $r-1$
vertices $x_1 \in V_1,\ldots,x_{r-1} \in V_{r-1}$ and element $s \in
S_i$, an edge with color $i$ and label $s$, which contains the
vertices $x_1,\ldots,x_{r-1}$ as well as vertex $y \in U_i$, where
\begin{equation}\label{simplevertex}
y=s + \sum^{r-1}_{j=1}a_{i,j}x_j\;,
\end{equation}
and the values $a_{i,j}$ were defined above. These edges will later
play the role of the edges $e_1,\ldots,e_{p-\ell}$ of $K$ defined
above. Note that these edges are defined irrespectively of the set
of equations $Mx=b$.

We now define the edges of $H$ which will ``simulate'' the linear
equations of $Mx=b$. For every $1 \leq i \leq \ell$, and for every
choice of an element $s \in S_{d_i}$, for every choice of
$r-1-|I_i|$ vertices $\{x_t \in V_t ~:~ t \in [r-1] \setminus I_i\}$
and for every choice of $|W_i| ~(=|I_i|)$ vertices $\{y_j \in U_j
~:~ j \in W_i\}$ we have an edge with color $d_i$ and label $s$,
which contains the vertices $\{x_t : t \in [r-1] \setminus I_i\}$
and $\{y_j~:~ j \in W_i\}$ as well as vertex $y \in U_{m_i}$, where
\begin{equation}\label{keyvertex}
y=b_i-M_{i,d_i} \cdot s - \sum_{j \in W_i}M_{i,j}\cdot y_j + \sum_{t
\in [r-1] \setminus I_i}x_t \cdot (a^{m_i}_t+\sum_{j \in W_i}a^j_{t}
\cdot M_{i,j})\;.
\end{equation}

Let us first note that as required by Lemma \ref{mainstep}, each
edge of $H$ has a color $i$ and is labeled by an element $s \in
S_i$. In fact, for each $1 \leq i \leq p$ and for each $s \in S_i$,
the hypergraph $H$ has $n^{r-1}$ edges that are colored $i$ and
labeled with $s$. We start with the following claim.

\begin{claim}\label{hsimple}
$H$ is a simple hypergraph, that is, it contains no parallel edges.
\end{claim}

\paragraph{Proof:}
Observe that edges of $H$ with different colors have a single vertex
from a different subset of $r$ of the sets
$V_1,\ldots,V_{r-1},U_1,\ldots,U_{p-\ell}$. Indeed, edges with color
$1 \leq i \leq p-\ell$ contain a vertex from each of the sets
$V_1,\ldots,V_{r-1}$ and another vertex from $U_i$, while an edge
with color $p-\ell+1 \leq d_i \leq p$ contains vertices from the
sets $\{V_t~:~ t \in [r-1]\setminus I_i\}$ as well as vertices from
some of the sets ${U_1,\ldots,U_{p-\ell}}$. Note that the sets
$I_1,\ldots,I_{\ell}$ are disjoint and non-empty, as none of the
sets $W_i$ is empty, a fact which (as noted previously) follows from
Claim \ref{notwo}. Observe that if $W_i$ was empty, then edges with
color $d_i$ would have had parallel edges with color $m_i$.

As for edges with the same color $1 \leq i \leq p-\ell$, recall that
they are defined in terms of a different combination of
$x_1,\ldots,x_{r-1} \in \mathbb{F}_n$ and $s \in S_i$. So if one edge
is defined in terms of $x_1,\ldots,x_{r-1} \in \mathbb{F}_n$ and $s
\in S_i$ and another using $x'_1,\ldots,x'_{r-1} \in \mathbb{F}_n$ and
$s' \in S_i$ then either (i) $x_j \neq x'_j$ for some $1 \leq j \leq
r-1$ in which case the edges have a different vertex in $V_j$ (ii)
$x_j=x'_j$ for all $1 \leq j \leq r-1$, implying that $s \neq s'$.
Therefore the edges have a different vertex in $U_i$ by the way we
chose the vertex in this set in (\ref{simplevertex}).

The case of edges with the same color $p-\ell+1 \leq d_i \leq p$ is
similar. Recall that such edges are defined in terms of a different
combination of $\{x_t ~:~ t \in  [r-1] \setminus I_i\}$, $\{y_j ~:~
j \in W_i\}$ and $s \in S_{d_i}$. So if one edge is defined in terms
of $\{x_t ~:~ t \in  [r-1] \setminus I_i\}$, $\{y_j ~:~ j \in W_i\}$
and $s \in S_{d_i}$ and another using $\{x'_t ~:~ t \in  [r-1]
\setminus I_i\}$, $\{y'_j ~:~ j \in W_i\}$ and $s' \in S_{d_i}$ then
either (i) $x_t \neq x'_t$ for some $t \in [r-1] \setminus I_i$ in
which case the edges have a different vertex in $V_t$ (ii) $y_j \neq
y'_j$ for some $j \in W_i$, in which case the edges have a different
vertex in $U_j$ (iii) $x_t = x'_t$ for all $t \in [r-1] \setminus
I_i$, and $y_j=y'_j$ for all $j \in W_i$,  implying that $s \neq s'$
and therefore the edges have a different vertex in $U_{m_i}$ by the
way we chose the vertex in this set in (\ref{keyvertex}) and from
the fact that $M_{i,d_i} \neq 0$. $\qed$

\bigskip

The above claim establishes the first property required by Definition \ref{DefHyperRep},
and we now turn to establish the second and third. Fix arbitrary elements $s_1 \in
S_1,\ldots,s_{p-\ell} \in S_{p-\ell}$. For every choice of $r-1$
(not necessarily distinct) elements $x_1,\ldots x_{r-1} \in
\mathbb{F}_n$, let $K_x$ be the set of vertices $x_1 \in V_1, \ldots,
x_{r-1} \in V_{r-1},y_1 \in U_1,\ldots,y_{p-\ell} \in U_{p-\ell}$,
where for every $1 \leq j \leq p-\ell$
\begin{equation}\label{defys}
y_j=s_j+\sum^{r-1}_{t=1}a_t^j \cdot x_t \;.
\end{equation}

We will need the following important claim regarding the vertices of
$K_x$. Getting back to the overview of the proof given in Subsection
\ref{subsecover}, this is where we extract one of the linear
equations $L_i$ (defined above) from a certain combination of edges
of a copy of $K$. We also note that the linear equation we
``initially'' obtain (see (\ref{importanteq})) includes also the
elements $x_i$, but the way we have constructed $H$ guarantees that
the $x_i$'s vanish and we eventually get a linear equation involving
only elements from the sets $S_i$. We will then use this claim to
show that $H$ contains many edge disjoint copies of $K$ when
$s_1,\ldots,s_{p-\ell}$ determine a solution to $Mx=b$, and in the
other direction, that $H$ cannot contain too many copies of $H$. For
what follows we remind that reader that for $1 \leq i \leq \ell$ we
have $p-\ell+1 \leq d_i \leq p$ and that for $i <i'$ we have $d_i
\neq d_{i'}$. Returning to the overview of the proof given in
Subsection \ref{subsecover}, we are now going to use the fact that
edges with colors $d_i$ and $m_i$ have a common vertex in $U_{m_i}$
in order to deduce the linear equation $L_i$.

\begin{claim}\label{edgeequation}
Let $1 \leq i \leq \ell$. Then the vertices $\{x_t ~:~ t \in
[r-1]\setminus I_i\}\cup \{y_j ~:~ j \in W_i\} \cup y_{m_i}$ span an
edge (of color $d_i$) if and only if there is an element $s_{d_i}
\in S_{d_i}$ such that $\{s_j~:~ j \in W_i\} \cup s_{m_i} \cup
s_{d_i}$ satisfy equation $L_i$ (defined in (\ref{rewriteeq})).
\end{claim}

\paragraph{Proof:}
$H$ contains an edge containing the vertices $\{x_t ~:~ t \in
[r-1]\setminus I_i\}\cup \{y_j ~:~ j \in W_i\} \cup y_{m_i}$ if and
only if (recall (\ref{keyvertex})) there is an $s_{d_i} \in S_{d_i}$ such that
\begin{equation}\label{importanteq}
y_{m_i}=b_i-M_{i,d_i} \cdot s_{d_i} - \sum_{j \in W_i}M_{i,j}\cdot y_j +
\sum_{t \in [r-1] \setminus I_i}x_t\cdot (a^{m_i}_t+\sum_{j \in
W_i}a^j_{t} \cdot M_{i,j})
\end{equation}
Using (\ref{defys}) this is equivalent to
requiring that
\begin{eqnarray*}
s_{m_i}+\sum^{r-1}_{t=1}a^{m_i}_{t} \cdot x_t
&=& b_i-M_{i,d_i} \cdot s_{d_i} - \sum_{j \in W_i}M_{i,j}\cdot (s_j+\sum^{r-1}_{t=1}a^j_{t} \cdot x_t)\\
&&~~+ \sum_{t \in [r-1] \setminus I_i} x_t \cdot (a^{m_i}_t+\sum_{j \in W_i}a^j_{t} \cdot M_{i,j}) \\
&=& b_i-M_{i,d_i} \cdot s_{d_i} - \sum_{j \in W_i}M_{i,j}\cdot s_j -  \sum^{r-1}_{t=1} x_t \cdot \left(\sum_{j \in W_i}a^j_{t} \cdot M_{i,j} \right)  \\
&&~~+ \sum_{t \in [r-1] \setminus I_i}x_t \cdot (a^{m_i}_t+\sum_{j \in W_i}a^j_{t} \cdot M_{i,j})\\
&=& b_i-M_{i,d_i} \cdot s_{d_i} - \sum_{j \in W_i}M_{i,j}\cdot s_j -
\sum_{t \in I_i} x_t \cdot \left(\sum_{j \in W_i}a^j_{t} \cdot
M_{i,j} \right)\\&&~~ + \sum_{t \in [r-1] \setminus I_i}x_t \cdot
a^{m_i}_t.
\end{eqnarray*}
Using (\ref{nicesum}) in the last row above, we can write the above requirement as
$$
s_{m_i}+\sum^{r-1}_{t=1}a^{m_i}_{t} \cdot x_t=
b_i-M_{i,d_i} \cdot s_{d_i} - \sum_{j \in W_i}M_{i,j}\cdot s_j
+\sum^{r-1}_{t=1}a^{m_i}_{t} \cdot x_t\;,
$$
or equivalently that
$$
s_{m_i}+M_{i,d_i}\cdot s_{d_i} +\sum_{j \in W_i}M_{i,j}\cdot
s_{j}=b_i\;,
$$
which is precisely equation $L_i$. $\qed$

\bigskip

For the next two claims, let us recall that we assume that the last
$\ell$ columns of $M$ form a diagonal matrix. Therefore, a solution
to $Mx=b$ is determined by the first $p-\ell$ elements of $x$.

\begin{claim}\label{copypersol}
Suppose $s_1,\ldots,s_{p-\ell}$ determine a solution
$s_1,\ldots,s_{p}$ to $Mx=b$. Then, any set $K_x$ (defined above)
spans a colored copy of $K$. In particular, for every solution
$s_1,\ldots,s_p$ to $Mx=b$, $H$ has $n^{r-1}$ colored copies of $K$,
in which the edge of color $i$ is colored with $s_i$.
\end{claim}

\paragraph{Proof:}
We claim that $K_x$ spans a colored copy of $K$, where for every $1
\leq i \leq r-1$ vertex $v_i$ of $K$ is mapped to vertex $x_i$ of
$H$, and for every $1 \leq j \leq p-\ell$ vertex $u_j$ of $K$ is
mapped to vertex $y_j$ of $H$. To see that the above is a valid
mapping of the colored edges of $K$ to colored edges of $H$, we
first note that the way we have defined $H$ in (\ref{simplevertex})
and the vertices $y_1,\ldots,y_{p-\ell}$ in (\ref{defys}),
guarantees that for every $1 \leq j \leq p-\ell$ we have an edge
with color $i$ which contains the vertices $x_1,\ldots,x_{r-1},y_j$.
This is actually true even if $s_1,\ldots,s_{p-\ell}$ do not
determine a solution.

As for edges with color $p-\ell+1 \leq d_i \leq p$, the fact that
the vertices $\{x_t ~:~ t \in [r-1]\setminus I_i\}\cup \{y_j ~:~ j
\in W_i\} \cup y_{m_i}$ span such an edge follows from Claim
\ref{edgeequation}, because we assume that $s_1,\ldots,s_{p-\ell}$
determine a solution to $Mx=b$, so for every $1 \leq i \leq \ell$
there exists an element $s_{d_i} \in S_{d_i}$ as required by Claim
\ref{edgeequation}. We thus conclude that
$x_1,\ldots,x_{r-1},y_1,\ldots,y_{p-\ell}$ span a colored copy of
$K$. Finally, note that by the way we have defined $H$, the edge of
$K_x$ which is colored $i$ is indeed labeled with the element $s_i
\in S_i$. $\qed$

\begin{claim}\label{copiesdisjoint}
If $s_1,\ldots,s_{p-\ell}$ determine a solution to $Mx=b$, then the
$n^{r-1}$ colored copies of $K$ spanned by the sets $K_x$ (defined
above) are edge disjoint.
\end{claim}

\paragraph{Proof:}
Let us consider two colored copies $K_x$ and $K_y$ for some $x \neq
y$ (Claim \ref{copypersol} guarantees that $K_x$ and $K_y$ indeed
span a colored copy of $K$). Clearly $K_x$ and $K_y$ cannot share
edges with color $1 \leq i \leq p-\ell$, because the vertices of
such edges within $V_1,\ldots,V_{r-1}$ are uniquely determined by
the coordinates of $x$ and $y$.

We now consider an edge of $K_x$ with color $d_i \in \{p-\ell+1,\ldots ,p\}$.
Let $j_1 < j_2 < \ldots < j_{|W_i|}$ be the elements of
$W_i$, and let $B_i$ be the matrix defined in Claim \ref{Bi}. 
Recall that $B_i$ satisfies the following\footnote{We remark
that when we have defined the matrices $B_i$ in Claim \ref{Bi} we did not
``impose'' the ordering of the rows that correspond to $W_i$ as we do here,
but this ordering, of course, does not affect the rank of $B_i$.}: (i) for $j \in [r-1] \setminus I_i$
we have $(B_i)_{j,j}=1$ and $(B_i)_{j,t}=0$ when $t \neq j$, and (ii) 
if $j \in I_i$ is the $g^{th}$ element of $I_i$, then the $j^{th}$ 
row of $B_i$ is the vector $a^{j_g}$ (where $j_g$ is the $g^{th}$ element of $W_i$). 
Let us also define an $r-1$ dimensional vector $c$ as follows: for every $j \in [r-1] \setminus I_i$ 
we have $c_j=0$, and for every $j \in I_i$, if $j$ is the $g^{th}$ element of $I_i$
then $c_j=s_{j_g}$. The key observation now is that the
vertices of the edge whose color is $d_i \in \{p-\ell+1, \ldots,
p\}$ within the $r-1$ sets $\{V_j ~:~ j \in [r-1] \setminus I_i\}
\cup \{U_j ~:~ j \in W_i\}$ are given by $B_ix+c$. More precisely,
for every $j \in [r-1] \setminus I_i$ the vertex of the edge of
color $d_i$ within $V_j$ is given by $(B_ix+c)_j$. Also,
for every $j_g \in W_i$, if $j \in I_i$ is the $g^{th}$ element of $I_i$,
then the vertex of this edge within $U_{j_g}$ is
given by $(B_ix+c)_j$. Claim \ref{Bi} asserts
that $B_i$ is non-singular, so we can conclude that the edges with
color $d_i$ of $K_x$ and $K_y$ can share at most $r-2$ of their
$r-1$ vertices within the sets $\{V_j ~:~ j \in [r-1]\setminus I_i\}
\cup \{U_j ~:~ j \in W_i\}$. So any pair of edges of color $d_i$ can
share at most $r-1$ vertices, and therefore $K_x$ and $K_y$ are edge
disjoint \footnote{We note that the way we have defined $H$ does not
(necessarily) guarantee that edges of the same color cannot share
$r-1$ vertices. That is, edges of color $i$ may share the vertex in
the set $U_{m_i}$ and $r-2$ of the $r-1$ vertices from the sets
$\{V_j ~:~ j \in [r-1]\setminus I_i\} \cup \{U_j ~:~ j \in W_i\}$.}.
$\qed$

\begin{claim}\label{copynum}
If $S_1,\ldots,S_p$ contain $T$ solutions to $Mx=b$ with $x_i \in
S_i$ then $H$ contains $Tn^{r-1}$ colored copies of $K$.
\end{claim}

\paragraph{Proof:} Recall that we assume that the last $\ell$ columns
of $M$ form a diagonal matrix. Therefore, the number of solutions
$T$ to $Mx=b$ is just the number of choices of $s_1 \in
S_1,\ldots,s_{p-\ell} \in S_{p-\ell}$ that can be extended to a
solution of $Mx=b$ by choosing appropriate values $s_{p-\ell+1} \in
S_{p-\ell+1},\ldots,s_p \in S_p$. Therefore, it is enough to show
that every colored copy of $K$ in $H$ is given by a choice of $r-1$
vertices $x_1 \in V_1,\ldots,x_{r-1}\in V_{r-1}$ and a choice of
$p-\ell$ elements $s_1 \in S_1,\ldots,s_{p-\ell} \in S_{p-\ell}$
that determine a solution to $Mx=b$.
So let us consider a colored copy of $K$ in $H$. This copy must
contain edges with the colors $1,\ldots,p-\ell$. By the way we have
defined $H$ this means that this copy must contain $r-1$ vertices
$x_1 \in V_1, \ldots, X_{r-1} \in V_{r-1}$ as well as $p-\ell$
vertices $y_1 \in U_1, \ldots, y_{p-\ell} \in U_{p-\ell}$. Furthermore,
for $1 \leq j \leq p-\ell$ we have
\begin{equation}\label{step1}
y_j = s_j + \sum^{r-1}_{t=1}a^j_t \cdot x_t
\end{equation}
for some choice of $s_j \in S_j$. So the vertex set of such a copy
is determined by the choice of $x_1,\ldots,x_{r-1}$ and
$s_1,\ldots,s_{p-\ell}$. Note that the set of vertices is just the
set $K_x$ defined before Claim \ref{edgeequation}, for
$x_1,\ldots,x_{r-1}$ and $s_1,\ldots,s_{p-\ell}$. Therefore, we can
apply Claim \ref{edgeequation} on this set of vertices.

So our goal now is to show that there are elements
$s_{p-\ell+1},\ldots,s_p$ which together with
$s_1,\ldots,s_{p-\ell}$ form a solution of $Mx=b$. Consider any $1
\leq i \leq \ell$. As the vertices at hand span a colored copy of
$K$ they must span an edge with color $d_i$. This edge
must\footnote{Because only vertices from this combination of $r$ of
the sets $V_1,\ldots,V_{r-1},U_1,\ldots,U_{p-\ell}$ spans an edge
with color $d_i$.} contain the vertices $\{x_t~:~ t \in
[r-1]\setminus I_i\} \cup \{y_j ~:~ j \in W_i\} \cup y_{m_i}$. But
by Claim \ref{edgeequation} if these vertices span an edge (of color
$d_i$) then there is an element $s_{d_i} \in S_{d_i}$ such that
$\{s_j~:~ j \in W_i\} \cup s_{m_i} \cup s_{d_i}$ satisfy equation
$L_i$. As this holds for every $1 \leq i \leq \ell$ we deduce that
$s_1,\ldots,s_p$ satisfy $Mx=b$. $\qed$

\bigskip

The proof of Lemma \ref{mainstep} now follows from Claims
\ref{hsimple}, \ref{copypersol}, \ref{copiesdisjoint} and \ref{copynum}.

\section{Concluding Remarks and Open Problems}\label{concluding}

\begin{itemize}

%\item Polynomial equations

\item Our removal lemma for sets of linear equations works over
any field. For the special case of a single linear equation,
Kr\'al', Serra and Vena \cite{Kral} (following Green \cite{Green})
proved a removal lemma over any group. It is natural to ask if a
similar removal lemma over groups, or even just abelian groups, also
holds for sets of linear equations.

\item Green \cite{Green} used the regularity lemma for groups
in order to resolve a conjecture
of Bergelson, Host, Kra and Ruzsa \cite{BHK}, which stated that every
$S \subseteq [n]$ of size $\delta n$ contains at least $(\delta^3-o(1))n$
3-term arithmetic progressions with a common difference. The analogous statement
for arithmetic progressions of length more than 4 was shown to be false in
\cite{BHK}. So the only case left open is whether any $S \subseteq [n]$
of size $\delta n$ contains at least $(\delta^4-o(1))n$ 4-term arithmetic progressions
with a common difference. Part of the motivation of Green for raising
Conjecture \ref{Greensconj} was that it may help in resolving the case of
the 4-term arithmetic progression. It seems very interesting to see if
Theorem \ref{maintheo} can indeed help in resolving this conjecture.

\item Our proof of the removal lemma for sets of linear equations
applies the hypergraph removal lemma. As a consequence, we get
extremely poor bounds relating $\epsilon$ and $\delta$. Roughly
speaking, the best current bounds for the graph removal lemma give
that $\delta(\epsilon)$ grows like Tower$(1/\epsilon)$, that is, a
tower of exponents of height $1/\epsilon$. For 3-uniform
hypergraphs, the bounds are given by iterating the Tower
function $1/\epsilon$ times, and so on. So on the one hand, the fact
that we are using hypergraphs with a large degree of uniformity
implies that the bounds we get are are extremely weak. On the other
hand, as even the graph removal lemma gives bounds which are too
weak for any reasonable application, this is not such a real issue
to be concerned about. It may still be interesting, however, to see
if one can prove Theorem \ref{maintheo} with a proof similar
to the one given in \cite{Kral} for the special case of a single equation.

\item Given the above discussion it it reasonable to ask for which
sets of equations $Mx=b$ one can get a polynomial dependence between
$\epsilon$ and $\delta$. This seems to be a challenging open problem
even for a single equation so let us focus on this case. For a
linear equation $L$, let $r_L(n)$ denote the size of the largest
subset of $n$ which contains no (non-trivial) solution to $L$.
Problems of this type were studied by Ruzsa \cite{R}. A simple
counting argument shows that if $r_L(n)=n^{1-c}$ for some positive
$c$, then $\delta(\epsilon)=O(1/\epsilon)^{1/c}$. However,
characterizing the equations with this property seems like a very
hard problem, see \cite{R}. Furthermore, we do not even know if all
the linear equations for which $r_L(n)=n^{1-o(1)}$ do not have a
polynomial dependence between $\epsilon$ and $\delta$. For example,
we do not know if such a dependence exists for the linear equation
$x_1+x_2=x_3$ (for which $r_L(n)=\Theta(n)$).

But for at least some of these linear equations, we can rule out
such a polynomial dependence as the following example shows.
Consider the linear equation $x_1+x_3=2x_2$, that is, the linear
equation which defines a 3-term arithmetic progression\footnote{The
argument can be extended to any linear equation in which one
variable is a convex combination of the others.}. We claim that for
this equation there is no polynomial relation between $\epsilon$ and
$\delta$. Fix an $\epsilon$ and let $n_0=n_0(\epsilon)$ be large
enough so that every $S \subseteq [n]$ of size $\epsilon n$ contains
a 3-term arithmetic progression. Roth's Theorem \cite{R} states that
such an $n$ exists. Therefore, for every $n \geq n_0$ and for every
$S \subseteq [n]$ of size $2\epsilon n$ we have to remove at least
$\epsilon n$ elements from $S$ in order to destroy all 3-term
arithmetic progressions. Let $m$ be the largest integer for which
$[m]$ contains a subset of size $4\epsilon m$, containing no 3-term
arithmetic progressions. The well known construction of Behrend
\cite{B} implies that $m \geq (1/\epsilon)^{c\log (1/\epsilon)}$ for
some absolute constant $c$. Let $X$ be one such subset of $[m]$. For
every $n \geq n_0$, let $S \subseteq [n]$ be the set of integers
with the property that in their base $2m$ representation, the least
significant element belongs to $X$. Then clearly
$|S|=n\cdot\frac{|X|}{2m}=2\epsilon n$ and so one should remove at
least $\epsilon n$ elements from $S$ to destroy all 3-term
arithmetic progressions. On the other hand if $x_1,x_2,x_3 \in S$
form a 3-term arithmetic progression then as $X \subseteq [m]$, so
do the least significant characters of $x_1,x_2,x_3$, because there
in no carry in the base $2m$ addition. But as these characters
belong to $X$ we get that they must be identical. Therefore, the
number of 3-term arithmetic progressions in $S$ is $|S|^3/m^2 \leq
\epsilon^{c\log 1/\epsilon }n^3$, implying that $\delta(\epsilon)
\leq \epsilon^{c\log 1/\epsilon}$.

\item The contrapositive version of our main result says that if one
should remove $\epsilon n$ elements from $S \subseteq [n]$ in order
to destroy all solutions of $Mx=b$ then $S$ contains
$f(\epsilon)n^{p-\ell}$ solutions to $Mx=b$. The ``analogous'' result
for graphs (or hypergraphs) is that if one should remove
$\epsilon n^2$ edges from a graph $G$ in order to destroy all the
copies of $H$ then $G$ contains $\delta(\epsilon)n^h$ copies of $H$ (where
$h$ is the number of vertices of $H$). The main result of \cite{AS}
is an ``infinite'' version of the removal lemma for graphs, which
states that if ${\cal H}$ is a (possibly infinite) set of graphs,
and if one should remove $\epsilon n^2$ edges from $G$ in order to
destroy all the copies of all the graphs $H \in {\cal H}$ then for
some $H \in {\cal H}$, whose size $h$ satisfies $h \leq h(\epsilon)$, $G$ contains
$\delta(\epsilon)n^h$ copies of $H$. It
seems natural to ask if there is a corresponding ``infinite" removal
lemma for sets of linear equations. More precisely, is it the case
that for every (possibly infinite) set ${\cal
M}=\{M_1x=b_1,M_2x=b_2,\ldots\}$ of sets of linear equations the
following holds: if one should remove $\epsilon n$ elements from $S
\subseteq [n]$ in order to destroy all the solutions to all the sets
of linear equations in ${\cal M}$, then for some set of linear
equations $Mx=b \in {\cal M}$, with $p \leq p(\epsilon)$ unknowns, $S$ contains
$\delta(\epsilon)n^{p-\ell}$ solutions to $Mx=b$.

\end{itemize}

\paragraph{Acknowledgements:} We would like to thank Vojta R\"odl, Benny
Sudakov and Terry Tao for helpful discussions related to this paper. I would
also like to thank Pablo Candela for his helpful comments on the paper.

\end{document}